\documentclass[12pt]{article}
\usepackage{amsfonts,amsmath,amsxtra}
\usepackage{latexsym}
\usepackage[matrix,arrow,curve]{xy}
\usepackage{amssymb}
\def\hybrid{\topmargin 0pt      \oddsidemargin 0pt
        \headheight 0pt \headsep 0pt
        \textwidth 16.5cm
        \textheight 23cm
        \marginparwidth 0.0in
        \parskip 5pt plus 1pt   \jot = 1.5ex}
\catcode`\@=11
\def\marginnote#1{}
\newcount\hour
\newcount\minute
\newtoks\amorpm
\hour=\time\divide\hour by60 \minute=\time{\multiply\hour by60
\global\advance\minute by-\hour}
\edef\standardtime{{\ifnum\hour<12 \global\amorpm={am}%
        \else\global\amorpm={pm}\advance\hour by-12 \fi
        \ifnum\hour=0 \hour=12 \fi
      \number\hour:\ifnum\minute<10 0\fi\number\minute\the\amorpm}}
\edef\militarytime{\number\hour:\ifnum\minute<10 0\fi\number\minute}

\def\draftlabel#1{{\@bsphack\if@filesw {\let\thepage\relax
   \xdef\@gtempa{\write\@auxout{\string
      \newlabel{#1}{{\@currentlabel}{\thepage}}}}}\@gtempa
   \if@nobreak \ifvmode\nobreak\fi\fi\fi\@esphack}
        \gdef\@eqnlabel{#1}}
\def\@eqnlabel{}
\def\@vacuum{}
\def\draftmarginnote#1{\marginpar{\raggedright\scriptsize\tt#1}}

\def\draft{\oddsidemargin -0.1truein
        \def\@oddfoot{\sl preliminary draft \hfil
        \rm\thepage\hfil\sl\today\quad\militarytime}
        \let\@evenfoot\@oddfoot \overfullrule 3pt
        \let\label=\draftlabel
        \let\marginnote=\draftmarginnote
\def\@eqnnum{{\rm (\theequation)}
\rlap{\kern\marginparsep\tt\@eqnlabel}%
\global\let\@eqnlabel\@vacuum}  }

\newfont{\Bbbb}{msbm7 scaled 1\@ptsize00}
\newcommand{\zs}{\raise-1pt\hbox{$\mbox{\Bbbb Z}$}}

\def\numberbysection{\@addtoreset{equation}{section}
        \def\theequation{\thesection.\arabic{equation}}}
\numberbysection

\renewcommand{\theequation}{\thesection.\arabic{equation}}

\def\titlepage{\@restonecolfalse\if@twocolumn\@restonecoltrue\onecolumn
     \else \newpage \fi \thispagestyle{empty}\c@page\z@
\def\thefootnote{\fnsymbol{footnote}} }
\def\endtitlepage{\if@restonecol\twocolumn \else  \fi
        \def\thefootnote{\arabic{footnote}}
        \setcounter{footnote}{0}}  
\relax
\hybrid
\makeatletter
\newdimen\normalarrayskip            
\newdimen\minarrayskip               
\normalarrayskip\baselineskip \minarrayskip\jot
\newif\ifold             \oldtrue            \def\new{\oldfalse}
\def\arraymode{\ifold\relax\else\displaystyle\fi}
\def\eqnumphantom{\phantom{(\theequation)}} 
\def\@arrayskip{\ifold\baselineskip\z@\lineskip\z@
     \else
     \baselineskip\minarrayskip\lineskip1\baselineskip\fi}
\def\@arrayclassz{\ifcase \@lastchclass \@acolampacol \or
\@ampacol \or \or \or \@addamp \or
   \@acolampacol \or \@firstampfalse \@acol \fi
\edef\@preamble{\@preamble
  \ifcase \@chnum
     \hfil$\relax\arraymode\@sharp$\hfil
     \or $\relax\arraymode\@sharp$\hfil
     \or \hfil$\relax\arraymode\@sharp$\fi}}
\def\@array[#1]#2{\setbox\@arstrutbox=\hbox{\vrule
     height\arraystretch \ht\strutbox
     depth\arraystretch \dp\strutbox
width\z@}\@mkpream{#2}\edef\@preamble{\halign \noexpand\@halignto
\bgroup \tabskip\z@ \@arstrut \@preamble \tabskip\z@ \cr}%
\let\@startpbox\@@startpbox \let\@endpbox\@@endpbox
  \if #1t\vtop \else \if#1b\vbox \else \vcenter \fi\fi
  \bgroup \let\par\relax
  \let\@sharp##\let\protect\relax
  \@arrayskip\@preamble}
\def\eqnarray{\stepcounter{equation}%
              \let\@currentlabel=\theequation
              \global\@eqnswtrue
              \global\@eqcnt\z@
              \tabskip\@centering              
              \let\\=\@eqncr
              $$%
            \halign to \displaywidth  \bgroup
             \eqnumphantom \@eqnsel
      \hskip\@centering                               
    $\displaystyle  \tabskip\z@ {##}$%
    &\global\@eqcnt\@ne \hskip 2\arraycolsep
         $ \displaystyle  \arraymode{##}$\hfil
    &\global\@eqcnt\tw@ \hskip 2\arraycolsep
         $\displaystyle\tabskip\z@{##}$\hfil
         \tabskip\@centering
    &{##}\tabskip\z@\cr}
\makeatother

\def\IC{\mathbb{C}}

\def\IR{\mathbb{R}}
\def\IZ{\mathbb{Z}}

\def\CH {\mathcal{H}}

\def\CU {\mathcal{U}}
\def\CV {\mathcal{V}}

\def\CZ {\mathcal{Z}}
\def\Fg{{\frak g}}

\def\fg{\mathfrak{g}}
\def\fh{\mathfrak{h}}
\def\fn{\mathfrak{n}}

\def\a {{\alpha}}

\def\g {{\gamma}}

\def\la{\lambda}

\def\ve{\varepsilon}

\def\e{\epsilon}

\def\fh{\mathfrak{h}}

\def\gl{\mathfrak{gl}}

\def\0{\bar{0}}
\def\1{\bar{1}}
\def\str{{\rm str}\,}

\def\pr {\partial}

\def\Id{{\rm Id}}

\def\End{{\rm End}}

\def\frak{\mathfrak}

\def\tr{{\rm tr}\,}

\def\<{\langle}
\def\>{\rangle}

\def\ad{{\rm ad}}

\def\osp{\mathfrak{osp}}

\newtheorem{de}{Definition}[section]
\newtheorem{prop}{Proposition}[section]           

\newtheorem{lem}{Lemma}[section]

\newtheorem{rem}{Remark}[section]

\newcommand\bqa{\begin{eqnarray}}
\newcommand\eqa{\end{eqnarray}}
\def\be{\begin{eqnarray}\new\begin{array}{cc}}
\def\ee{\end{array}\end{eqnarray}}
\def\beq{\begin{equation}}
\def\eeq{\end{equation}}
\def\bse{\begin{subequations}}                
\def\ese{\end{subequations}}
\def\bp{\begin{pmatrix}}
\def\ep{\end{pmatrix}}

\newcommand{\proof}{\noindent {\it Proof} }

\newcounter{pac}[section]

\newcounter{pacc}[subsection]

0

\begin{document}

\title{\bf On quantum $\osp(1|2\ell)$-Toda chain
  \footnote{Talk given by the second author at the "Polivanov-90"
    conference,
    16-17 December 2020, Steklov Mathematical Institute of Russian Academy of Science.}}
\author{A.A. Gerasimov, D.R. Lebedev and S.V. Oblezin}
\date{}
\maketitle

\renewcommand{\abstractname}{}

\begin{abstract}
\noindent {\bf Abstract}. The orthosymplectic superalgebra
$\mathfrak{osp}(1|\,2\ell)$ is the  closest analog of standard Lie
algebras in the world of super Lie algebras. We demonstrate that the
corresponding $\osp(1|\,2\ell)$-Toda chain turns out to be an
instance of a $BC_\ell$-Toda chain. The underlying reason for this
relation is discussed.
\end{abstract}



\vspace{5 mm}

\section{Introduction}

Representation theory is an essential tool in finding explicit
solutions of known quantum integrable systems as well as
construction of new ones. An important  class of  finite-dimensional
quantum integrable systems allowing representation theory
interpretation is provided by Toda chains. It is  known that
integrable Toda chains  are classified by a class of root systems
that include root systems of finite dimensional Lie algebras as well
as their affine counterparts. For the  Toda chains associated with
the root systems of finite dimensional Lie algebras the
corresponding integrable systems can be solved explicitly by
representation theory methods \cite{K1}, \cite{GW} (see \cite{STS}
for a review). Precisely, the eigenfunctions of the quantum
Hamiltonians are given by special matrix elements of principal
series representations of the totally split real form of the
corresponding Lie group. The resulting functions should  be
considered as generalized Whittaker functions associated with the
corresponding finite-dimensional  Lie algebras \cite{K1}, \cite{Ha}.
These functions allow quite explicit integral representations (see
e.g. \cite{GLO1}).

The class of integrable Toda chains is however a bit larger than the
class of finite/affine Lie algebras and includes in particular
non-reduce  root systems $BC_\ell$ combining $B_\ell$ and $C_\ell$
root systems. The corresponding $BC_\ell$-Toda system is an
important element of the web of Toda type theories connected by
various intertwining relations \cite{GLO2}. Although $BC_\ell$ root
system fits naturally in the classification of finite Lie algebra
root systems the problem of  construction of the Lie algebra type
object corresponding to the non-reduced root $BC_\ell$ systems seems
not  yet obtained a  satisfactory resolution. However, one should
recall that $BC_\ell$ root systems appear in the Cartan
classification of symmetric spaces \cite{H}, \cite{L}. Still
$BC_\ell$-Toda chain can be solved via representation theory methods
using  a generalization of $C_\ell$-Whittaker functions (see e.g.
\cite{J} for $\ell=1$ and a remark in \cite{RS}, relevant to
$BC_\ell$ classical Toda system). This unfortunately  does not
elucidate the question of the
 interpretation of $BC_\ell$-Toda eigenfunctions as standard  Whittaker
functions for some group-like object. One should add that the
integrability of the quantum $BC_\ell$-Toda chain for generic coupling
constants was  proven independently in \cite{S} using Yangian
representation theory (aka quantum inverse scattering methods).
This however also does not clarify the question of existence of
a group-like structure behind $BC_\ell$ root systems.

In this note we consider quantum Toda chains associated with the
super Lie algebras $\osp(1|2\ell)$.  This series of super Lie
algebras  occupy a special place in the world of  super
Lie algebras. In particular, it is the only instance  of  simple
super Lie algebras for which the corresponding category of
finite-dimensional representations  is semi-simple
and thus allows direct analogs of  the standard constructions of
representation theory of semisimple Lie algebras \cite{Kac1}.
In connection with this fact one should mention that
$\osp(1|2\ell)$ is the unique super Lie algebra
with finitely-generated center of its universal enveloping algebra.
The special properties of $\osp(1|2\ell)$ makes it natural to consider
the associated quantum integrable  systems.

In this note we demonstrate  that $\osp(1|2\ell)$-Toda chain
may be also considered  as a Toda chain associated with the $BC_\ell$
root system.  This allows us to solve $BC_\ell$-Toda chain
by standard representation theory methods i.e. by identifying
the corresponding eigenfunctions with $\osp(1|2\ell)$-Whittaker functions.

The underlying reason  for the appearance of $BC_\ell$
root structure in  $\osp(1|2\ell)$-Toda chain becomes clear by
comparing $BC_\ell$ root data with that of the super Lie algebra
$\osp(1|2\ell)$. Actually the only difference is the opposite parity
of the maximal commutative subalgebra eigenspaces in the Cartan
decomposition  corresponding to short roots
of non-reduced $BC_\ell$ root system. This difference however does not
affect the expressions for quantum Hamiltonians of the
corresponding Toda chain.

Exposition of the paper goes as follows. In Section 2 we provide
basic facts on the structure of the orthosymplect super Lie algebra
$\osp(1|2\ell)$. In Section 3 we construct $\osp(1|2\ell)$-Whittaker
functions  associated with representations of the
super Lie algebra $\osp(1|2\ell)$ and demonstrate that these functions are
eigenfunctions of  the quadratic quantum Hamiltonian of $BC_\ell$-Toda
chain for special values of the coupling constants. Finally,
in Section 4 we discuss the structure of root system of
$\osp(1|2\ell)$ versus  $BC_\ell$ root system and provide an
explanation of the  apparent  identification of
quadratic Hamiltonian of $\osp(1|2\ell)$-Toda chain with
that of  $BC_\ell$-Toda chain.

{\it  Acknowledgments:} The research of the second (D.R.L.) and
third (S.V.O.) authors  was supported by RSF grant 16-11-10075.

\section{Basic facts on the super Lie algebra $\osp(1|2\ell)$}

We start with the basic definition of a Lie superalgebra structure
and then we describe explicitly the algebra  $\osp(1|2\ell)$ in
detail. This is a standard material that can be found in standard
sources on super algebras e.g. \cite{Kac1}, \cite{Kac2}.

The notion of Lie superalgebra is a direct generalization of the
notion of Lie algebra  to the category of vector superspaces. Vector
superspace  $V=V_{\0}\oplus V_{\1}$ is a $\IZ_2$-graded vector space
with the parity $p$ taking values $0$ and $1$ on $V_{\0}$ and
$V_{\1}$ respectively. The tensor product structure is given by
twisting of the standard tensor product structure in the category of
vector spaces \be v\otimes w=(-1)^{p(v)\cdot p(w)}\,w\otimes
v,\qquad v\in V, \quad w\in W, \ee for $v$ and $w$ are homogeneous
elements with respect to the $\IZ_2$-grading.

\begin{de} The structure of super Lie algebra on super vector space
  $\fg=\fg_{\0}\oplus\fg_{\1}$ is given by
 a bilinear operation $[\cdot,\,\cdot]$, called the bracket, so
that for any homogeneous elements $X,\,Y,\,Z\in \fg$ the following hold:
 \be
  p\bigl([X,\,Y]\bigr)\,=\,p(X)\,+\,p(Y)\,,
 \ee
 \be
  \,[X,\,Y]\,=\,-(-1)^{p(X)\cdot p(Y)}[Y,\,X]\,,
 \ee
 \be
  \,[X,\,[Y,\,Z]](-1)^{p(X)\cdot p(Z)}\,
  +\,[Z,\,[X,\,Y]](-1)^{p(Y)\cdot p(Z)}\\
  +\,[Y,\,[Z,\,X]](-1)^{p(X)\cdot p(Y)}\,=\,0\,.
 \ee
\end{de}

We will  be interested in a special instance  of super Lie algebras,
the ortho-symplectic super Lie algebra  $\osp(1|2\ell)$. To define
this algebra let us first introduce the super Lie algebra
$\gl(1|2\ell)$.

\begin{de}\label{GL1} The super Lie algebra $\gl(1|\,2\ell)$ is generated by
 \be
  E_{0,\,i}\,,\quad E_{i,\,0}\,,\quad
  p(E_{0,\,i})\,=\,p(E_{i,\,0})\,=\,1\,,\qquad0\leq i\leq2\ell\,,\\
  \text{and}\qquad E_{kl}\,,\qquad p(E_{kl})\,=\,0\,,\qquad1\leq
  k,\,l\leq2\ell\,,
 \ee
subjected to the following relations:
 \be\label{glbracket}
  \bigl[E_{ij},\,E_{kl}\bigr]\,
  =\,\delta_{jk}E_{il}\,-(-1)^{p(i)p(l)}\,\delta_{il}E_{kj}\,,\qquad
  0\leq i,\,j\leq2\ell\,,\quad
  0\leq k,\,l\leq2\ell\,.
 \ee
\end{de}

The super Lie agebra $\gl(1|2\ell)$ may be identified with the super
Lie algebra structure $\bigl(\End(V),\,[\cdot,\,\cdot]\bigr)$
 on the space ${\rm End}(V)$  of endomorphisms of the superspace
\be
V\,=\,\IR^{1|2\ell}\,=\,V_{\0}\,\oplus\,V_{\1}\,,\qquad
 V_{\0}\,=\,\IR^{0|2\ell}\,,\qquad V_{\1}\,=\,\IR^{1|0}\,,
 \ee
 in the following way.
 Any zero parity linear endomorphism $A\in\End(V)$ is given by the  matrix of the
following shape:
 \be\label{shape}
  A\,  =\,\left(\begin{array}{cc}
   A_{11} & A_{12}\\A_{21} & A_{22}
  \end{array}\right)\,,\quad
  \begin{array}{cc}
  A_{11}:\,V_{\1}\,\longrightarrow\,V_{\1}\,, &
  A_{12}:\,V_{\0}\,\longrightarrow\,V_{\1}\,,\\
  A_{21}:\,V_{\1}\,\longrightarrow\,V_{\0}\,, &
  A_{22}:\,V_{\0}\,\longrightarrow\,V_{\0}\,,
  \end{array}
 \ee
where  entries of blocks $A_{11},\,A_{22}$ are even while  the
entries of $A_{12}\,,A_{21}$ are odd so that
 \be
  \End(V)_{\0}\,=\,\Big\{\Big(\begin{smallmatrix}
  A_{11}&&0\\&&\\0&&A_{22}\end{smallmatrix}\Big)\Big\}\,,\qquad
  \End(V)_{\1}\,=\,\Big\{\Big(\begin{smallmatrix}
  0&&A_{12}\\&&\\A_{21}&&0\end{smallmatrix}\Big)\Big\}\,.
 \ee
The super brackets on ${\rm End}(V)$ are defined on homogeneous elements $X,\,Y\in\End(V)$
as follows
 \be\label{bracket}
  [X,\,Y]\,=\,X\circ Y\,-\,(-1)^{p(X)\cdot p(Y)}Y\circ X\,.
  \ee
The description of $\fg(1|2\ell)$ given in Definition
\ref{GL1} is then obtained via fixing  a bases in $V$
\be\label{basis}
\{\ve_0,\,\ve_1,\,\ldots,\,\ve_{2\ell}\}\subset\IR^{1|2\ell}\,,\qquad
p(\ve_0)=1\,,\quad p(\ve_k)=0,\quad1\leq k\leq2\ell\,.
 \ee
 The generators $E_{ij}$ are identified with the elementary matrices
 in ${\rm End}(V)$  with the only non-zero elements  in the $i$-th row
and the  $j$-th column.

\begin{de} The super transposition of a matrix $A\in\End(V)$ is
defined by
 \be\label{sTransp}
  A^{\top}\,
  =\,\left(\begin{array}{cc}
   A_{11} & A_{12}\\A_{21} & A_{22}
  \end{array}\right)^{\top}\,
  =\,\left(\begin{array}{cc}
   A_{11}^t & -A_{21}^t\\A_{12}^t & A_{22}^t
  \end{array}\right)\,,
 \ee
where $X^t$ if the standard transposition of a matrix $X$.
\end{de}
\begin{lem} Super transposition \eqref{sTransp} possesses the
following properties:
 \be
  (A\,v)^t\,=\,v^tA^{\top}\,,\qquad A\in\End(V)\,,\quad v\in V\,,
 \ee
 \be
  (A\cdot B)^{\top}\,=\,B^{\top}\cdot A^{\top}\,,
 \ee
 \be
  (A^{\top})^{\top}\,=\Pi\,A\,\Pi^{-1}\,,
 \ee
where $\Pi$ is the parity operator with the matrix
$\Big(\begin{smallmatrix}-1&&0\\&&\\0&&\Id_{2\ell}\end{smallmatrix}\Big)$.
\end{lem}
\proof: Given $v\in V$ let us write it down in the basis
$\{\ve_0,\,\ve_1,\ldots,\ve_{2\ell}\}$:
 \be
  v\,=\,\xi\ve_0\,+\,\sum_{i=1}^{2\ell}v_i\ve_i\,,
 \ee
with odd Grassmann coordinate $\xi$, and even coordinates $v_i$.
Then we have
 \be
  (A\,v)^t_i\,=\,a_{i0}\xi\,+\,a_{i1}v_1\,+\ldots+\,a_{i,\,2\ell}v_{2\ell}\,,
 \ee
and on the other hand,
 \be
  (v^tA^{\top})_0\,
  =\,\xi a_{00}\,+\,v_1a_{0,1}\,+\ldots+\,v_na_{0,n}\,,\\
  (v^tA^{\top})_k\,=\,-\xi
  a_{k0}\,+\,v_1a_{k,1}\,+\ldots+\,v_na_{k,n}\,,\quad1\leq k\leq2\ell\,.
 \ee
Taking into account that
 \be
  \xi a_{00}\,=\,a_{00}\xi\,,\qquad-\xi
  a_{k,0}\,=\,a_{k,0}\xi\,,\quad1\leq k\leq2\ell\,,
 \ee
we deduce the first assertion. The second assertion can be verified
by straightforward computation. The third assertion follows from the
definition: on the one hand, we have
 \be
  (A^{\top})^{\top}\,
  =\,\left(\begin{array}{cc}
   A_{11}^t & -A_{21}^t\\A_{12}^t & A_{22}^t
  \end{array}\right)^{\top}\,
  =\,\left(\begin{array}{cc}
   A_{11} & -A_{12}\\-A_{12} & A_{22}
  \end{array}\right)\,;
 \ee
on the other hand, in the standard basis \eqref{basis} the matrix of
parity operator reads
 \be
  \Pi\,=\,\Big(\begin{smallmatrix}
  -1 && 0\\&&\\0&&\Id_{2\ell}
  \end{smallmatrix}\Big)\,,
 \ee
so the assertion easily follows. $\Box$

Now  $\osp(1|2\ell)$ may be defined as a subalgebra of the  general linear
superalgebra  $\gl(1|2\ell)$. Introduce the following involution:
 \be
  \theta\,:\quad\gl(1|2\ell)\,\longrightarrow\,\gl(1|2\ell)\,,\qquad
  X\,\longmapsto\,X^{\theta}\,:=\,-JX^{\top}J^{-1}\,,
  \ee
where
 \be
  J\,
  =\,\Big(\begin{smallmatrix}
  1 & 0 & 0\\0 & 0 & -\Id_{\ell}\\
  0 & \Id_{\ell} & 0
  \end{smallmatrix}\Big)\,\in\,\End(V)_{\0}\,.
 \ee

\begin{de} The orthosymplectic super Lie algebra $\osp(1|2\ell)$ is defined as
the $\theta$-invariant subalgebra of $\mathfrak{gl}(1|2\ell)$:
 \be\label{ospNrep}
  \osp(1|\,2\ell)\,
  =\Big\{X\,\in\,\gl(1|2\ell)\,:\quad X^{\theta}\,=\,X\Big\}\\
  =\,\Big\{X\,
  =\,\Big(\begin{smallmatrix}
   0 & x & y\\
   y^t & A & B\\
   -x^t & C & -A^t
  \end{smallmatrix}\Big)\,:\quad
  B^t\,=\,B\,,\quad
  C^t\,=\,C\Big\}\,\subset\,\gl(1|\,2\ell)\,.
 \ee
\end{de}

According to the classification of simple super Lie algebras
\cite{Kac1} one associates the root system $B_{0,\ell}$ to the super
Lie algebra $\osp(1|2\ell)$. Let
$\{\e_1,\ldots,\,\e_{\ell}\}\subset\IR^{\ell}$ be an orthogonal
basis in $\IR^{\ell}$ with respect to the scalar product $(\,,\,)$.
Then simple root system ${}^s\Delta^+(B_{0,\ell})$ of type
$B_{0,\,\ell}$ consists of even simple positive roots
${}^s\Delta^+_{\bar{0}}$ and odd simple positive roots
${}^s\Delta^+_{\bar{1}}$:
 \be\label{OPSroot}
  {}^s\Delta^+_{\bar{0}}(B_{0,\,\ell})\,=\,
  \bigl\{\a_k\,=\,\e_{\ell+1-k}-\e_{\ell+2-k}\,,\quad1<k\leq\ell\bigr\}\,\,,\\
  {}^s\Delta^+_{\bar{1}}(B_{0,\,\ell})\,=\,\bigl\{\a_1=\e_{\ell}\bigr\}\,,
 \ee
indexed by $I=\{1,\,\ldots,\,\ell\}$. The simple co-roots
$\a_i^{\vee},\,i\in I$ are defined in a standard way:
$$
 \a_i^{\vee}\,:=\,\frac{2\a_i}{(\a_i,\a_i)}\,,\quad i\in I\,.
$$
Note that the set $\Delta^+(B_{0,\ell})$ of positive roots contains the
sub-system of even positive roots of $C_\ell$ root system with the
corresponding set of simple roots:
 \be
  {}^s\Delta^+(C_{\ell})\,=\,\bigl\{2\a_1=2\e_{\ell}\,,\quad
  \a_k\,=\,\e_{\ell+1-k}-\e_{\ell+2-k}\,,\,\,1<k\leq\ell\bigr\}\,
  \subset\,\Delta^+(B_{0,\ell})\,.
 \ee

The Cartan matrix $A=\|A_{ij}\|$ associated with the simple
 root system \eqref{OPSroot} is defined by the standard
 formula
 \be
  A_{ij}\,=\,\frac{2(\alpha_i,\alpha_j)}{(\alpha_i,\alpha_i)}\,,\qquad i,j\in I\,.
 \ee
Thus the Cartan matrix of ${}^s\Delta^+(B_{0,\ell})$ coincides with
the standard $B_{\ell}$-type Cartan matrix
 \be\label{OSPcar}
  A\,
  =\,\left(\begin{array}{c|cccc}
   2&-2&0& \ldots &0\\
  \hline
   -1&2&-1&\ddots&\vdots\\
   0&\ddots&\ddots&\ddots&0\\
  \vdots&\ddots&-1&2&-1\\
   0&\ldots&0&-1&2
  \end{array}\right).
\ee
The Cartan decomposition for $\osp(1|2\ell)$ reads
 \be\label{ospNCartan}
  \osp(1|2\ell)(\IC)\,=\,\bigoplus_{i\in I}\IC
  h_i\,\oplus\,\bigoplus_{\a\in\Delta^+_{\bar{0}}}\bigl(\IC
  X_{\a}\,\oplus\,\IC X_{-\a}\bigr)\,\oplus\,\bigoplus_{\beta\in\Delta^+_{\bar{1}}}\bigl(\IC
  X_{\beta}\,\oplus\,\IC X_{-\beta}\bigr)\,,\\
  \Delta^+_{\bar{0}}\,=\,\bigl\{2\e_i\,;\qquad\e_i\pm\e_j\,,\quad
  i<j\,,\quad i,\,j\in I\bigr\}\,,\qquad
  \Delta^+_{\bar{1}}\,=\,\bigl\{\e_i\,,\quad i\in I\bigr\}\,.
 \ee
and the Cartan-Weyl relations are the following:
 \be\label{osprel}
  \bigl[X_{\e_i},\,X_{\e_j}\bigr]\,=\,(1+\delta_{ij})X_{\e_i+\e_j}\,,\qquad
  \bigl[X_{-\e_i},\,X_{-\e_j}\bigr]\,=\,-(1+\delta_{ij})X_{-\e_i-\e_j}\,,\\
  \bigl[X_{\e_i},\,X_{-\e_i}\bigr]\,=\,a_{ii}\,,\qquad i\in I\,;\\
  \bigl[X_{\e_i-\e_j},\,X_{\e_j}\bigr]\,=\,X_{\e_i}\,,\qquad
  \bigl[X_{\e_i-\e_j},\,X_{-\e_i}\bigr]\,=\,-X_{-\e_j}\,,\\
  \bigl[X_{\e_i},\,X_{-\e_i-\e_j}\bigr]\,=\,X_{-\e_j}\,,\qquad
  \bigl[X_{-\e_i},\,X_{\e_i+\e_j}\bigr]\,=\,X_{\e_j}\,,\qquad i<j\,;\\
  \bigl[X_{\a},\,X_{-\a}\bigr]\,=\,h_{\a^{\vee}}\,
  =\,\sum_{i\in I}\<\a^{\vee},\,\e_i\>a_{ii}\,,\\
  \bigl[h_{\a^{\vee}},\,X_{\g}\bigr]\,=\,\a^{\vee}(\g)X_{\g}\,,\qquad
  \a,\,\g\in\Delta^+\,.
 \ee
The Serre relations on $X_{\a_i},\,\a_i\in{}^s\Delta^+(B_{0,\ell})$
have the following form:
 \be\label{Serre}
  \ad_{X_{\a_1}}^2(X_{\a_1})\,=\,0\,,\qquad
  \ad_{X_{-\a_1}}^2(X_{-\a_1})\,=\,0\,,\\
  \ad_{X_{\a_i}}^{1-a_{ij}}(X_{\a_j})\,=\,0\,,\quad
  \ad_{X_{-\a_i}}^{1-a_{ij}}(X_{-\a_j})\,=\,0\,,\qquad
  i,j\in I\,\,.
 \ee

The Cartan-Weyl generators $X_\alpha$ may be represented via matrix
embedding \eqref{ospNrep} of
 $\osp(1|2\ell)$ as follows:
 \be\label{ospNgen1}
  X_{\e_i}\,=\,E_{i,\,0}\,+\,E_{0,\,\ell+i}\,,\qquad
  X_{-\e_i}\,=\,E_{0,\,i}\,-\,E_{\ell+i,\,0}\,;
 \ee
 \be\label{ospNgen0}
  X_{\e_i-\e_j}\,=\,E_{ij}-E_{2\ell+1-i,\,2\ell+1-j}\,,\\
  X_{-\e_i+\e_j}\,=\,E_{ji}-E_{2\ell+1-j,\,2\ell+1-i}\,,\\
  X_{\e_i+\e_j}\,=\,E_{i,\,\ell+j}+E_{j,\,\ell+i}\,,\quad
  X_{-\e_i-\e_j}\,=\,E_{\ell+j,\,i}+E_{\ell+i,\,j}\,,\quad
  i<j\,,\\
  X_{2\e_i}\,=\,E_{i,\,\ell+i}\,,\qquad
  X_{-2\e_i}\,=\,E_{\ell+i,\,i}\,,\qquad i\in I\,.
 \ee
The Cartan subalgebra $\fh\subset\osp(1|\,2\ell)$ is spanned by
 \be\label{ospNgen2}
  h_i\,=\,E_{ii}\,-\,E_{\ell+i,\,\ell+i}\,,\qquad i\in I\,.
 \ee

For a class of super Lie algebras $\fg$ allowing
non-degenerate invariant pairing $(\,|\,)$ there is a canonical
construction of the  quadratic Casimir element
$C_2\in\CZ(\CU(\fg))$ of the center of the universal enveloping
algebra $\CU(\fg)$ (see e.g.
\cite{Kac1}). Let us chose a  pair  $\{u_i,\,i\in I\}$,
$\{u^i,\,i\in I\}$ of  dual bases
in the Cartan subalgebra $\fh\subset\fg$, and let
$\{X_{\a},\,X^{\a},\,\a\in\Delta^+\}$ be the Cartan-Weyl generators
normalized by $(X^{\a}|\,X_{\a})=1$. Then the  quadratic Casimir element
$C_2$  allows for the following presentation:
 \be\label{casimir}
  C_2\,
  =\,\sum_{i\in I}u^iu_i\,
  +\,\sum_{\a\in\Delta^+}\bigl((-1)^{p(\a)}X_{\a}X^{\a}\,+\,X^{\a}X_{\a}\bigr)\,.
 \ee
To define  the quadratic Casimir element for $\osp(1|2\ell)$ we shall
first introduce an invariant non-degenerate invariant pairing.
The super Lie algebra $\osp(1|2\ell)$ allows the  invariant scalar
product  defined as follows
 \be\label{ISP2}
  (X|Y)\,:=\,\frac{1}{2}\str\bigl(\rho_t(X)\circ
  \rho_t(Y)\bigr)\,,\qquad X,\,Y\in \osp(1|2\ell)\,,
 \ee
where $\rho_t: \osp(1|2\ell)\to \End(\IC^{1|2\ell})$ is the
  tautological representation of $\osp(1|2\ell)$ in $\IC^{1|2\ell}$.
The supertrace of $A\in\End(\IC^{1|2\ell})$ of the shape
\eqref{shape} is given by
 \be
  \str(A)\,=\,\str\left(\begin{array}{cc}
   A_{11} & A_{12}\\A_{21} & A_{22}
  \end{array}\right)\,
  =\,-\tr(A_{11})\,+\,\tr(A_{22})\,.
 \ee
The explicit form of the invariant scalar product \eqref{ISP2}
may be directly derived using the matrix representation
\eqref{ospNrep}.

\begin{lem} The invariant scalar product \eqref{ISP2} on super Lie
  algebra  $\osp(1|2\ell)$ is as follows
 \be\label{ISP}
  (h_i|h_i)=1\,,\quad i\in I\,;\\
  (X_\alpha|X_{-\alpha})=\frac{2}{(\alpha,\alpha)}\,,\quad
  \alpha\in \Delta^+_{\bar{0}}\,;
  \qquad  (X_\beta|X_{-\beta})=1\,,
  \quad \beta\in \Delta^+_{\bar{1}}\,,
 \ee
with the rest of the products being zero.
\end{lem}

\proof: Validity of \eqref{ISP} may be checked directly. Thus  for
example we have using \eqref{ospNgen2}
 \be
  (h_i|h_j)=\frac{1}{2}\str\,\Big(\rho_t(E_{ii}-E_{\ell+i,\ell+i})
  \rho_t(E_{jj}-E_{\ell+j,\ell+j})\Big)\,
  =\,\delta_{ij}\,.
 \ee
Similarly, using \eqref{ospNgen1} we obtain
 \be
  (X_{\epsilon_i}|X_{-\epsilon_j})\,
  =\,\frac{1}{2}\str\Big(\rho_t(E_{i,0}+E_{0,\ell+i})
  \rho_t(E_{0,j}-E_{\ell+j,0})\Big)\,
  =\,\delta_{ij}\,.
 \ee
Similarly one might check the expressions for remaining products. $\Box$

\begin{prop} The following expression provides  the quadratic
Casimir element \eqref{casimir} for $\osp(1|2\ell)$:
 \be\label{Gcasimir}
  C_2\,
   =\,\sum_{i\in I}\Big(h_{i}^2\,
   -\,X_{\e_i}X_{-\e_i}\,+\,X_{-\e_i}X_{\e_i}\Big)\\
   +\,\sum_{\a\in\Delta^+_{\bar{0}}}\frac{(\a,\,\a)}{2}\bigl(X_{\a}X_{-\a}\,
   +\,X_{-\a}X_{\a}\bigr)\,.
  \ee
\end{prop}
\proof:  Using the expressions \eqref{ISP} for the invariant pairing
we have
 \be\label{ISP1}
  h^i\,=\,h_i\,,\quad i\in I\,;\\
  X^{\alpha}\,=\,\frac{(\alpha,\alpha)}{2}\,X_{-\alpha}\,,\quad
  \alpha\in \Delta^+_{\bar{0}}\,;\qquad
  X^{\beta}\,=\,X_{-\beta}\,,\quad \beta\in \Delta^+_{\bar{1}}\,.
 \ee
Substituting  \eqref{ISP1} into \eqref{casimir} we arrive at
\eqref{Gcasimir}, and complete the proof. $\Box$

In the following it will more convenient to use
another set of notations for generators which is adapted to the
matrix form \eqref{ospNrep}
 \be
  y_i\,=X_{\e_i}\,,\qquad   x_i\,=\,X_{-\e_i}\,\,;\\
  a_{ii}\,=\,h_i\,,
  \quad b_{ii}\,=\,X_{2\e_i}\,,\quad  c_{ii}\,=\,X_{-2\e_i}\,,\qquad i\in I\,;\\
  a_{ij}\,=\,X_{\e_i-\e_j}\,\qquad a_{ji}\,=\,X_{-\e_i+\e_j}\,,\\
  b_{ij}\,=\,X_{\e_i+\e_j}\,,\quad c_{ij}\,=\,X_{-\e_i-\e_j},\quad
  i<j\,\,.
 \ee
In addition to \eqref{osprel} the even part of $\osp(1|2\ell)$
satisfies the following relations:
 \be\label{spNrels}
  [b_{ij},b_{kl}]=0, \qquad
  [c_{ij},c_{kl}]=0, \qquad
  [b_{ij},c_{kl}]=\delta_{jk}a_{il},\\
  \,[a_{ii},b_{kl}]=(\delta_{ik}-\delta_{il})b_{kl}, \qquad
  [a_{ii},c_{kl}]=-(\delta_{ik}-\delta_{il})b_{kl}\,.
 \ee
Using these notations the quadratic Casimir element \eqref{Gcasimir}
may be written as follows:
 \be\label{Gcasimir1}
   C_2\,
  =\,\sum_{i=1}^{\ell}(a_{ii}^2+x_iy_i-y_ix_i)\,
  +\,2(c_{ii}b_{ii}+b_{ii}c_{ii})\\
  +\,\sum_{i<j}(a_{ij}a_{ji}+a_{ji}a_{ij})\,+\,(b_{ij}c_{ij}+c_{ij}b_{ij})\,.
  \ee
From now on we will consider the real form $\osp(1|2\ell)(\IR)$ of
the orthosymplectic super Lie algebra, such that the generators
$a_{ii},\,b_{ii},\,c_{ii},\,i\in I$,
$a_{ij},\,a_{ji},\,b_{ij},\,c_{ij},\,i<j$ as well as $x_i$ and $y_i$
are defined to be real.

\section{The $\osp(1|\,2\ell)$-Whittaker function}

In this section we construct the Whittaker function associated with
the super Lie algebra $\osp(1|2\ell)$. There is a classical approach
to the construction of  Whittaker functions associated with
semisimple Lie algebras \cite{J}, \cite{K1}, \cite{K2}, \cite{H}.
Below we give a modified version of this construction due to
Kazhdan-Kostant (see \cite{E}).

Given a super Lie algebra $\Fg$, let $\CU(\Fg)$ be the corresponding
universal enveloping algebra and let $\CZ(\CU(\fg))\subset\CU(\Fg)$
be its center. A  $\CU(\fg)$-module $\CV$ admits an infinitesimal
character $\zeta$ if there is a homomorphism $\zeta:
\CZ(\CU(\fg))\rightarrow \IC$ such that $zv=\zeta(z)v$ for all
$z\!\in\!\CZ(\CU(\fg))$ and $v\in\CV$. Given a character $\chi$ of a
nilpotent super Lie subalgebra $\fn\subset \fg$,
 \be
  \chi\,:\quad\fn\longrightarrow\,\IC^{1|1}\,;
  \ee
we define a Whittaker vector $\psi\in\CV$  by the following
relations:
 \be
  X\cdot\psi\,=\,\chi(X)\,\psi\,,\qquad\forall
  X\in\fn\subset\fg\,.
 \ee
The  Whittaker vector $\psi\in\CV$ is  called cyclic, if it
generates $\CV$: $\CU(\Fg)\,\cdot\psi=\CV$. A $\CU(\Fg)$-module
$\CV$ is called a Whittaker module if it contains a cyclic Whittaker
vector. A pair of $\CU(\Fg)$-modules  $\CV$ and $\CV'$ is called
dual if there exists a non-degenerate  pairing
 \be
  \<\,,\,\>\,:\quad\CV\times\CV'\,\longrightarrow\,\IC^{1|1}\,,
 \ee
which is $\IC$-antilinear in  the first
 variable and $\IC$-linear in the second one, and  such that
 \be\label{hermitean}
  \<X\cdot v',\,v\>\,=\,-(-1)^{p(v')\cdot p(X)}\<v'\,,X\cdot
  v\>\,,\qquad v\in\CV,\,\, v'\in\CV'\,,\quad X\in\fg\,.
 \ee

Now we restrict ourselves to the case of the orthosymplectic super
Lie algebra $\osp(1|2\ell)$. Let $\CV_\la$ be a $\CU(\osp(1|\,2\ell))$-module
with an infinitesimal central character allowing
a vector $v_{\la}\in\CV_{\la}$ defined by (see the notations
of \eqref{ospNCartan}, \eqref{osprel}):
 \be\label{hwvec}
  h_{\a^{\vee}}\cdot v_{\la}\,=\,\a^{\vee}(\la)v_{\la}\,,\quad
  X_{\a}\cdot v_{\la}\,=\,0\,,\qquad\forall
  X_{\a}\in\fn_+\subset\osp(1|2\ell)\,,\quad\a\in\Delta^+\,,
 \ee
 where $\la$ is an element of the dual to the Cartan subalgebra
 $\mathfrak{h}\subset \osp(1|2\ell)$.
Value of the quadratic Casimir element $C_2$ on $\CV_{\la}$ is
uniquely determined by \eqref{hwvec}. Indeed let us re-write the
Casimir element $C_2$ from \eqref{Gcasimir} as follows:
 \be\label{Gcasimir11}
  C_2\,=\,\sum_{i\in I}\Big(a_{ii}^2\,
   -\,a_{ii}\,+\,2X_{-\e_i}X_{\e_i}\Big)\,
   +\,\sum_{\a\in\Delta^+_{\bar{0}}}\frac{(\a,\,\a)}{2}\Big(h_{\a^{\vee}}\,
   +\,2X_{-\a}X_{\a}\Big)\\
  =\,\sum_{i\in I}\Big(a_{ii}^2\,+\,2\rho(\e_i)a_{ii}\,
   +\,2X_{-\e_i}X_{\e_i}\Big)\,
   +\,\sum_{\a\in\Delta^+_{\bar{0}}}(\a,\,\a)\,X_{-\a}X_{\a}\,,
 \ee
where
 \be\label{rho1}
  \rho(q)\,
  =\,\frac{1}{2}\Big(\sum_{\a\in\Delta^+_{\bar{0}}}\a(q)\,
  -\,\sum_{\beta\in\Delta^+_{\bar{1}}}\beta(q)\Big)\,.
 \ee
Thus using \eqref{hwvec} $C_2$ takes the following value on
$v_{\la}\in\CV_{\la}$:
 \be\label{CasimirValue}
  C_2(v_{\la})\,=\,\sum_{i\in
  I}\bigl(a_{ii}^2\,+\,2\rho(\e_i)a_{ii}\bigr)(v_{\la})\,
  =\,(\la,\la+2\rho)\,v_{\la}\,.
 \ee

In the following we consider those $\CV_\lambda$ that allow a
structure of the  Whittaker modules and also allow integration of
the action of the Cartan subalgebra $\fh\subset\osp(1|2\ell)$ to the
action of the corresponding maximal torus $H$. Precisely let
$\CV_\lambda$ and $\CV'_\lambda$ be  a dual pair of Whittaker
modules  cyclically generated by
 Whittaker vectors $\psi_R\in\CV_{\lambda}$ and $\psi_L\in\CV_{\lambda}'$.
Explicitly, the Whittaker vectors $\psi_R\in\CV_{\lambda}$ and
$\psi_L\in\CV'_\lambda$ are defined by the following conditions
 \be\label{whittchar}
  X_{\a}\cdot\psi_R\,=\,\chi_R(X_{\a})\,\psi_R\,,\qquad
  X_{-\a}\cdot\psi_L\,=\,\chi_L(X_{-\a})\,\psi_L\,,\qquad
  \forall
  \a\in\Delta^+\,,
 \ee
where $\chi_R:\,\fn_+\to\IC^{1|1}$ and $\chi_L:\,\fn_-\to\IC^{1|1}$
are the characters of the opposite nilpotent super Lie subalgebras
$\fn_{\pm}\subset\osp(1|2\ell)$:
 \be
  \chi_R\,:\quad\fn_+\,=\,\Big(\bigoplus_{\a\in\Delta^+_{\bar{0}}}\IC
  X_{\a}\,\oplus\,\bigoplus_{\beta\in\Delta^+_{\bar{1}}}\IC
  X_{\beta}\Big)\,\longrightarrow\,\IC^{1|1}\,,\\
  \chi_L\,:\quad\fn_-\,=\,\Big(\bigoplus_{\a\in\Delta^+_{\bar{0}}}\IC
  X_{-\a}\,\oplus\,\bigoplus_{\beta\in\Delta^+_{\bar{1}}}\IC
  X_{-\beta}\Big)\,\longrightarrow\,\IC^{1|1}\,.
  \ee

\begin{lem}\label{Nchar}
\begin{itemize}
\item[(i)] The function $\chi_R:\,\fn_+\to\IC^{1|1}$ defined by
 \be\label{whittvecR}
  \chi_R\bigl(X_{\e_{\ell}}\bigr)\,=\,\imath^{3/2}\xi_{\a_1}^+\,\in\,
  \imath^{3/2}\IR^{1|0}\,,\quad
  \chi_R\bigl(X_{2\e_{\ell}}\bigr)\,=\,\imath(\xi_{\a_1}^+)^2\,\in\,\imath\IR^{0|1}\,,\\
  \chi_R\bigl(X_{\e_{\ell+1-k}-\e_{\ell+2-k}}\bigr)\,=\,
  \imath\xi_{\a_k}^+\,\in\,\imath\IR^{0|1}\,,
  \quad1<k\leq\ell\,,\\
  \chi_R\bigl(X_{\e_k}\bigr)\,=\,\chi_R\bigl(X_{\a}\bigr)\,=\,0\,,\qquad
  1\leq k<\ell\,,\quad\a\in\Delta^+_{\bar{0}}\setminus{}^s\Delta^+_{\bar{0}}\,,
 \ee
is a character of the super Lie subalgebra
$\fn_+\subset\osp(1|\,2\ell)$.

\item[(ii)] Similarly, the function $\psi_L:\,\fn_-\to\IC^{1|1}$ defined by
 \be\label{whittvecL}
  \chi_L\bigl(X_{-\e_{\ell}}\bigr)\,=\,\imath^{3/2}\xi_{\a_1}^-\,\in\,
  \imath^{3/2}\IR^{1|0}\,,\qquad
  \chi_L\bigl(X_{-\e_{2\ell}}\bigr)\,=\,\imath(\xi_{\a_1}^-)^2\,\in\,\imath\IR^{0|1}\,,\\
  \chi_L\bigl(X_{-\e_{\ell+1-k}+\e_{\ell+2-k}}\bigr)
  \,=\,\imath\xi_{\a_k}^-\,\in\,\imath\IR^{0|1}\,,\quad1<k\leq\ell\,,\\
  \chi_L\bigl(X_{-\e_k}\bigr)\,=\,\chi_L\bigl(X_{-\a}\bigr)\,=\,0\,,\qquad
  1\leq k<\ell\,,\quad\a\in\Delta^+_{\bar{0}}\setminus{}^s\Delta^+_{\bar{0}}\,,
 \ee
is a character of the super  subalgebra Lie
$\fn_-\subset\osp(1|\,2\ell)$.
\end{itemize}
\end{lem}
\proof: We provide the proof  in the case of $\chi_R$ while  the case of
$\chi_L$ can be treated in a similar way. Let us verify that
\eqref{whittvecR} defines a character $\chi_R$ of the super
subalgebra $\fn_+\subset\osp(1|2\ell)$ by checking the compatibility
of \eqref{whittvecR} with the appropriate Cartan-Weyl relations
\eqref{osprel}:
 \be\label{osprelN}
  [X_{\e_i},\,X_{\e_i}]\,=\,2X_{\e_i}^2\,=\,2X_{2\e_i}\,,\quad
  [X_{\e_i},\,X_{\e_j}]\,=\,X_{\e_i+\e_j}\,,\qquad i,j\in I\,,\\
  \bigl[X_{\e_i-\e_{i+1}},X_{\e_{i+1}}\bigr]=X_{\e_i},\quad
  \bigl[X_{\e_i-\e_{i+1}},[X_{\e_i-\e_{i+1}},X_{2\e_{i+1}}]\bigr]
  =2X_{2\e_i},\quad1\leq
  i<\ell\,,\\
  \bigl[X_{\e_i-\e_j},\,X_{2\e_j}\bigr]\,=\,X_{\e_i+\e_j}\,,\qquad i<j\,,\\
  \bigl[X_{\e_i-\e_j},\,X_{\e_j-\e_k}\bigr]\,=\,X_{\e_i-\e_k}\,,\qquad
  i<j<k\,,
 \ee
and with the Serre relations \eqref{Serre}:
 \be\label{SerreN}
  \ad_{X_{\e_{\ell}}}^2(X_{\e_{\ell}})\,=\,0\,,\quad
  \ad_{X_{\e_{\ell}}}^3(X_{\e_{\ell-1}-\e_{\ell}})\,
  =\,\ad_{X_{\e_{\ell-1}-\e_{\ell}}}^2(X_{\e_{\ell}})\,=\,0\,,\\
  \ad_{X_{2\e_{\ell}}}^2(X_{\e_{\ell-1}-\e_{\ell}})\,
  =\,\ad_{X_{\e_{\ell-1}-\e_{\ell}}}^3(X_{2\e_{\ell}})\,=\,0\,,\\
  \ad_{X_{\e_{i-1}-\e_i}}^2(X_{\e_i-\e_{i+1}})\,
  =\,\ad_{X_{\e_i-\e_{i+1}}}^2(X_{\e_{i-1}-\e_i})\,=\,0\,,\quad1<i<\ell\,.
 \ee

From the defining relations \eqref{whittvecR} we see that $\chi_R$
takes non-zero values only on the simple root generators
$X_{\e_{\ell}}$, $X_{\e_k-\e_{k+1}},\,1\leq k<\ell$ and on the
special non-simple root generator $X_{2\e_{\ell}}$. The latter
follows from the first relation from \eqref{osprelN} for $i=\ell$:
 \be
  [X_{\e_{\ell}},\,X_{\e_{\ell}}]\,=\,2X_{\e_{\ell}}^2\,
  =\,2X_{2\e_{\ell}}\,.
 \ee
Indeed, given $X_{\e_{\ell}}\cdot\psi_R=\imath^{3/2}\,\xi_{\a_1}^+\psi_R$ one
readily deduces
 \be
  2X_{\e_{\ell}}^2\cdot\psi_R\,=\,2X_{\e_{\ell}}\cdot(X_{\e_{\ell}}\cdot\psi_R)\,
  =\,2X_{\e_{\ell}}\cdot(\imath^{3/2}\xi_{\a_1}^+\psi_R)\,
  =\,-2\imath^{3/2}\xi_{\a_1}^+(X_{\e_{\ell}}\cdot\psi_R)\\
  =\,2\imath(\xi_{\a_1}^+)^2\psi_R\,,
 \ee
which matches with
$2X_{\e_{2\ell}}\cdot\psi_R\,=\,2\imath(\xi_{\a_1}^+)^2\psi_R$.
Similarly, the first relation from \eqref{osprelN} for $1\leq
i<\ell$ yields
 \be
  \chi_R(X_{2\e_i})\,=\,\frac{1}{2}\chi_R(X_{\e_i})^2\,
  =\,0\,,\qquad1\leq i<\ell\,,
 \ee
and the other relation in the first line of \eqref{whittvecR}
entails
 \be
  \chi_R(X_{\e_i+\e_j})\,=\,\chi_R(X_{\e_i}X_{\e_j}+X_{\e_j}X_{\e_i})\,
  =\,0\,,\qquad i<j\,.
 \ee
The Serre relations imply that
$\dim\fn_+=|\Delta^+|=\ell^2+\ell=\ell(\ell+1)$. Thus the rest of
the defining relations \eqref{whittvecR} are provided by the fact
that given $X_{\a}\in\fn_+$ and $X_{\beta},X_{\g}\in\fn_+$, such
that $\a=\beta+\g$ and not both $\beta,\g$ are odd, we have
 \be
  \chi_R(X_{\a})\,
  =\,\chi_R\bigl(X_{\beta}X_{\g}-X_{\g}X_{\beta}\bigr)\,=\,0\,.
 \ee
Namely, the last line of \eqref{osprelN} for each $i<j<k$ with
$1<i-k<\ell$ implies
 \be
  \chi_R(X_{\e_i-\e_k})\,
  =\,\chi_R\bigl(X_{\e_i-\e_j}X_{\e_j-\e_k}
  -X_{\e_j-\e_k}X_{\e_i-\e_j}\bigr)\,
  =\,0\,.
 \ee
Then the second and the third lines of \eqref{osprelN} for each
$1\leq i<k\leq\ell$ we have
 \be
  \chi_R(X_{\e_i})\,
  =\,\chi_R\bigl(X_{\e_i-\e_k}X_{\e_k}-X_{\e_k}X_{\e_i-\e_k}\bigr)\,=\,0\,.
 \ee
Finally, we check that the remaining relations
 \be
  \bigl[X_{\e_i-\e_j},\,X_{2\e_j}\bigr]\,=\,X_{\e_i+\e_j}\,,\qquad
  [X_{\e_i-\e_j},\,X_{2\e_j}]\,=\,X_{\e_i+\e_j}\,,\\
  \ad_{X_{\e_i-\e_j}}^2(X_{2\e_j})\,
  =\,\bigl[X_{\e_i-\e_j},[X_{\e_i-\e_j},X_{2\e_j}]\bigr]\,
  =\,X_{2\e_i}\,,\qquad i<j\,,
 \ee
are consistent with the defining relations \eqref{whittvecR}. This
completes our proof. $\Box$

\begin{rem} Our choice of the characters \eqref{whittvecL}, \eqref{whittvecR}
in Lemma  \ref{Nchar} is
  compatible with the notion of a unitary operators in the case of
  super Hilbert spaces (see \cite{DM}).  Note however that in our case we do not require
the Hilbert space structure but only an invariant pairing.
\end{rem}

\begin{de} Let $\CV_\lambda$ and $\CV'_\lambda$ be a dual pair of
  cyclic   Whittaker modules with the action of the Casimir element given
  by \eqref{CasimirValue}.  The $\osp(1|2\ell)$-Whittaker function is defined by
 \be\label{Gwhittaker}
  \Psi_{\la}(e^q)\,
  =\,e^{-\rho(q)}\<\psi_L\,,\,e^{-h_q} \cdot \psi_R\>\,,\qquad
  h_q\,=\,\sum_{i\in I}q_ia_{ii}\,,
 \ee
where $\rho$ is the half-sum of positive even roots minus the half-sum
of positive odd roots, given by \eqref{rho1}.
\end{de}

\begin{prop} The $\osp(1|\,2\ell)$-Whittaker function \eqref{Gwhittaker} is
a solution to the following eigenvalue problem:
 \be\label{eigenvalue}
 \CH_2^{\osp(1|2\ell)}\cdot\Psi_{\la}(e^q)
 \,=\,-(\lambda+\rho)^2\,\Psi_{\la}(e^q)\,,
 \ee
 \be\label{BCham}
 \begin{array}{c}
  \CH_2^{\osp(1|2\ell)}\,
   =-\,\sum_{i\in I}\frac{\pr^2}{\pr q_i^2}\,
   +\,2\!\!\sum_{\a_i\in{}^s\!\Delta^+}\xi_{\a_i}^-\xi_{\a_i}^+\,e^{\a_i(q)}\,
  +\,4(\xi_{\a_1}^-\xi_{\a_1}^+)^2\,e^{2\a_1(q)}\,,
 \end{array}
 \ee
where $\rho$ is given by \eqref{rho1}, and
${}^s\!\Delta^+={}^s\!\Delta^+(B_{0,\ell})$ is defined in
\eqref{OPSroot}.
\end{prop}

\proof: On the one hand, by our construction we read from
\eqref{CasimirValue}:
 \be\label{CasimirValue1}
  \<\psi_L\,,e^{-h_q}\,C_2\,\psi_R\>\,
  =\,(\lambda,\lambda+2\rho)\,\<\psi_L\,,e^{-h_q}\,C_2\,\psi_R\>\,.
 \ee
On the other hand, the action of the Casimir element
$C_2\in\CZ(\CU(\osp(1|2\ell)))$  is equivalent to action on
\eqref{Gwhittaker} of a  certain second-order differential operator.
Namely, from \eqref{Gcasimir1} we take
 \be
  C_2\,
  =\,\sum_{i\in I}\Big(a_{ii}^2\,+\,2\rho(\e_i)a_{ii}\,
   +\,2X_{-\e_i}X_{\e_i}\Big)\,
   +\,\sum_{\a\in\Delta^+_{\bar{0}}}(\a,\,\a)\,X_{-\a}X_{\a}\,,
 \ee
and substituting this into \eqref{CasimirValue1} we obtain:
 \be
 \sum_{i\in I}
 \bigl\<\psi_L\,,e^{h_q}\bigl(a_{ii}^2\,+\,2\rho(\e_i)a_{ii}\bigr)\,\psi_R\bigr\>\\
  =\,\sum_{i\in I}\Big\{\frac{\pr^2}{\pr q_i^2}\,
  -\,2\rho(\e_i)\frac{\pr}{\pr
  q_i}\Big\}\<\psi_L\,,e^{-h_q}\cdot\psi_R\>\,.
 \ee
Taking into account the defining equations \eqref{whittvecL},
\eqref{whittvecR}  and
 the hermitian property \eqref{hermitean} of $\<\,,\,\>$ we find out
 \be
  2\sum_{i\in I}\bigl\<\psi_L\,,e^{-h_q}X_{-\e_i}X_{\e_i}\,\psi_R\bigr\>\\
  =\,-2\sum_{i\in I}e^{q_i}
  (-1)^{p(X_{-\e_i})\,p(\psi_L)}\bigl\<X_{-\e_i}\psi_L\,\,,e^{-h_q}\,X_{\e_i}\,
  \psi_R\bigr\>\,\\
  =\,-2(-1)^{p(X_{-\e_{\ell}})\,p(\psi_L)}e^{q_{\ell}}
  \bigl\<\imath^{3/2}\xi_{\a_1}^-\psi_L\,,e^{-h_q}\imath^{3/2}\xi_{\a_1}^+\psi_R\bigr\>\,\\
  =\,-2(-1)^{p(X_{-\e_{\ell}})\,p(\psi_L)}(-1)^{p(\xi_{\a_1}^+)\cdot p(\psi_L)}
  \imath^{-3/2}\imath^{3/2}\xi_{\a_1}^-\xi_{\a_1}^+
  e^{\a_1(q)}\bigl\<\psi_L\,,e^{-h_q}\cdot\psi_R\bigr\>\,\\
  =\,-2\xi_{\a_1}^-\xi_{\a_1}^+e^{\alpha_1(q)}\,\bigl\<\psi_L\,,e^{-h_q}\cdot
    \psi_R\bigr\>\,.
 \ee
Here we use the fact that $\<\,,\,\>$ is $\IC$-antilinear in the
first  variable and it is $\IC$-linear in  the second variable. In a similar way we
 derive
 \be
  \sum_{\a\in\Delta^+_{\bar{0}}}(\a,\,\a)\bigl\<\psi_L\,,e^{-h_q}X_{-\a}X_{\a}
  \psi_R\bigr\>\\
  =\,-\sum_{\a\in\Delta^+_{\bar{0}}}(\a,\,\a)\,e^{\a(q)}
  \bigl\<X_{-\a}\,\psi_L\,,
 e^{h_q}X_{\a}\psi_R\bigr\>\\
 =\,-2\sum_{i=2}^{\ell}\imath^{-1}\imath\xi_{\a_i}^-\xi_{\a_i}^+e^{\a_i(q)}
 \bigl\<\psi_L\,,e^{-h_q}\psi_R\bigr\>\\
  -\,4\imath^{-1}\imath(\xi_{\a_1}^-\xi_{\a_1}^+)^2e^{2\a_1(q)}
  \bigl\<\psi_L\,,e^{-h_q}\psi_R\bigr\>\,.
 \ee
Collecting the contributions above we obtain the following:
 \be
  \<\psi_L\,,e^{-h_q}\,C_2\,\psi_R\>\,
  =\,\Big\{\sum_{i\in I}\Big(\frac{\pr^2}{\pr q_i^2}\,
  -\,2\rho(\e_i)\frac{\pr}{\pr q_i}\Big)\\
  -\,4(\xi_{\a_1}^-\xi_{\a_1}^+)^2e^{2\a_1(q)}\,
  -\,2\sum_{i=1}^{\ell}\xi_{\a_i}^-\xi_{\a_i}^+\,e^{\a_i(q)}
  \Big\}
  \<\psi_L\,,e^{-h_q}\cdot\psi_R\>\,.
 \ee
Now we observe that
 \be
  e^{-\rho(q)}\frac{\pr}{\pr q_i}e^{\rho(q)}\,
  =\,\frac{\pr}{\pr q_i}\,+\,\rho(q)'_{q_i}\,
  =\,\frac{\pr}{\pr q_i}\,+\,\rho(\e_i)\,,\\
  e^{-\rho(q)}\frac{\pr^2}{\pr q_i^2}e^{\rho(q)}\,
  =\,\frac{\pr^2}{\pr q_i^2}\,+\,2\rho(q)'_{q_i}\frac{\pr}{\pr q_i}\,
  +\,\bigl(\rho(q)'_{q_i}\bigr)^2\\
  =\,\frac{\pr^2}{\pr q_i^2}\,+\,2\rho(\e_i)\frac{\pr}{\pr q_i}\,
  +\,\rho(\e_i)^2\,,
 \ee
hence we deduce the following:
 \be
  \sum_{i\in I}e^{-\rho(q)}\Big\{\frac{\pr^2}{\pr q_i^2}\,
  -\,2\rho(\e_i)\frac{\pr}{\pr q_i}\Big\}e^{\rho(q)}\\
  =\,\sum_{i\in I}\Big\{\frac{\pr^2}{\pr q_i^2}\,+\,2\rho(\e_i)\frac{\pr}{\pr q_i}\,
  +\,\rho(\e_i)^2\,
  -\,2\rho(\e_i)\Big(\frac{\pr}{\pr q_i}\,+\,\rho(\e_i)\Big)\Big\}\\
  =\,\sum_{i\in I}\frac{\pr^2}{\pr q_i^2}\,-\,\rho^2\,.
 \ee
Finally, we collect all the contributions and substitute them into
\eqref{CasimirValue1} to deduce the following:
 \be
  \Big\{\sum_{i\in I}\frac{\pr^2}{\pr q_i^2}\,
  -\,4(\xi_{\a_1}^-\xi_{\a_1}^+)^2e^{2\a_1(q)}\,
  -\,2\sum_{\a_i\in{}^s\!\Delta^+}\xi_{\a_i}^-\xi_{\a_i}^+\,e^{\a_i(q)}
  -\,\rho^2\Big\}\cdot\Psi_{\la}(e^q)\\
  =\,(\la,\la+2\rho)\,\Psi_{\la}(e^q)\,,
 \ee
where ${}^s\!\Delta^+={}^s\!\Delta^+(B_{0,\ell})$. This easily
entails the assertion \eqref{eigenvalue}. $\Box$

\begin{rem} In the special case $\lambda=\imath\mu-\rho$,  the eigenvalue equation
\eqref{eigenvalue} reads
 \be\label{BCham1}
  \CH_2^{\osp(1|\,2\ell)}\,\Psi_{\la}(e^q)\,
  =\,\mu^2\Psi_{\la}(e^q)\,,\\
  \CH_2^{\osp(1|\,2\ell)}\,
  =\,-\sum_{i\in I}\frac{\pr^2}{\pr q_i^2}\,
   +\,2\!\!\sum_{\a_i\in{}^s\Delta^+}\xi_{\a_i}^-\xi_{\a_i}^+\,e^{\a_i(q)}\,
  +\,4(\xi_{\a_1}^-\xi_{\a_1}^+)^2\,e^{2\a_1(q)}\,.
 \ee
\end{rem}

Let us introduce the corresponding couplings:
 \be
  g_i^2\,=\,\xi^-_{\a_i}\,\xi^+_{\a_i}\,, \qquad i\in I\,.
\ee

\begin{lem}\label{INDEP} The $\osp(1|2\ell)$-Whittaker function \eqref{Gwhittaker}
  depends on $\xi^{\pm}_{\a_i}$ via $g_{i}^2$, $i\in I$.
\end{lem}

\proof: The $\osp(1|2\ell)$-Whittaker function \eqref{Gwhittaker}
\be \Psi_{\la}(e^q|\xi^{\pm}_{\a_i})\,=\,
e^{-\rho(q)}\<\psi_L\,,\,e^{-h_q}\cdot\psi_R\>\,,\qquad
  h_q\,=\,\sum_{i\in I}q_ia_{ii}\,,
 \ee
satisfies  the following obvious relation: given
$Q=\exp\bigl\{\sum\limits_{i=1}^\ell \theta_i h_i\bigr\}\in H$
 \be
  \<\psi_L\,,\,Qe^{-h_q}Q^{-1}\cdot\psi_R\>\,
  =\,\<\psi_L\,,\,e^{-h_q}\cdot\psi_R\>\,.
 \ee
The adjoint action of $Q$ on the left and right Whittaker vectors
$\psi_R,\,\psi_L$ \eqref{whittvecL}, \eqref{whittvecR} changes them,
so that the eigenvalues $\xi_{\a_i}^{\pm}$ of the corresponding
$\mathfrak{n}_{\pm}$-characters are changed as follows:
 \be
  \xi_{\a_1}^{\pm}\longrightarrow \xi_{\a_1}^{\pm}
  e^{\pm\theta_{\ell}}\,,\qquad
  \xi_{\a_i}^{\pm}\longrightarrow \xi_{\a_i}^{\pm}
  e^{\pm(\theta_{\ell+1-i}-\theta_{\ell+2-i})}\,,\quad1<i\leq\ell\,.
 \ee
The invariance of the Whittaker function under this transformation
implies  that the $\osp(1|2\ell)$ -Whittaker function  depends on
$\xi_{\a_i}^{\pm}$ only via quadratic combinations
$\xi_{\a_i}^+\xi_{\a_i}^-$. $\Box$

\begin{lem} Let us consider a specialization of the
  $\osp(1|2\ell)$-Toda chain
  by taking arbitrary special values $g_i^2=\kappa_i^2\in \IR$  of
  the couplings.  Then by a linear change of variables $q_i$ one can bring the quadratic
Hamiltonian
  \be
 \begin{array}{c}
  \CH_2^{\osp(1|\,2\ell)}\,
   =\,-\sum_{i\in I}\frac{\pr^2}{\pr q_i^2}\,
   +\,2\sum_{i=2}^{\ell} \kappa_i^2\,e^{\a_i(q)}\,
  +\,2\kappa_1^2\,e^{\a_1(q)}\,
  +\,4\kappa_1^4\,e^{2\a_1(q)}\,,
 \end{array}
 \ee
 to the following canonical form:
 \be\label{BCcan}
\CH_2^{\osp(1|\,2\ell)}\,
  =\,-\sum_{i\in I}\frac{\pr^2}{\pr q_i^2}\,
  +\,\sum_{i=2}^{\ell}e^{\a_i(q)}\,
  +\,e^{\a_1(q)}\,+\,e^{2\a_1(q)}.
 \ee
\end{lem}

\proof: Indeed it is easy to check that the following transformation
of variables,
 \be
  q_{\ell}\,\longmapsto\,q_{\ell}\,-\,\ln2\,-\,\ln \kappa_1^2\,,\\
  q_k\,\longmapsto\,q_k\,-\,(\ell+1-k)\ln2\,-\,\ln \kappa^2_{k}\,
  +\,\ln \kappa^2_{\ell-k}\,,\qquad1\leq k<\ell\,,
  \ee
  applied to  \eqref{BCham}  gives \eqref{BCcan}.  $\Box$

\section{On $\osp(1|2\ell)$ as a  Lie algebra of type $BC_\ell$}

Let us recall the construction of the Toda chain associated with the
general root system. Let $\Delta$ be a rank $\ell$  root system
realized as a set of vectors in $V=\IR^{\ell}$. Chose an orthogonal
basis $\{\e_i,\,i\in I\}$ in $V$ and the dual basis $\{\e^i,\,i\in
I\}$ in $V^*$, both indexed by $I=\{1,\ldots,\ell\}$. Then elements
$q\in V^*$ allow decomposition $q=\sum\limits_{i=1}^\ell q_i\e^i$.
Let ${}^s\!\Delta^+$ be as a set of  simple positive roots in
$\Delta$. The quadratic quantum Hamiltonian of the Toda chain
associated with the root system $\Delta$ is given by
 \be\label{qHam}
  \CH_2^{\Delta^+}\,
  =\,-\sum_{i\in I}\frac{\pr^2}{\pr q_i^2}+ \sum_{\alpha\in{}^s\!\Delta^+}\,g^2_\alpha
  \,e^{\alpha(q)}
 \ee
with the  coupling constants $g_\alpha^2$. Note that the Hamiltonian
depends only on the structure of simple positive roots
${}^s\!\Delta^+$.

Now let us specialize this expression to the case of $BC_\ell$-root
system,   the unique non-reduced root system  satisfying basic
axioms of root systems of finite-dimensional Lie algebras (for a
description of $BC_\ell$ root system see e.g. \cite{H}, \cite{L}).
The set of simple positive roots of the $BC_{\ell}$ root system is
given by
 \be\label{BC2}
  {}^s\Delta^+(BC_{\ell})\,
  =\,\bigl\{2\e_{\ell};\,\, \e_{\ell}\,,\quad \e_{i}-\e_{i+1}\,,\,\, 1\leq
  i<\ell\bigr\}\,.
 \ee
The Cartan matrix $A=\|A_{ij}\|$ associated with the  set of simple
positive roots is defined via standard formula
 \be\label{CAM}
  A_{ij}=\frac{2(\alpha_i,\alpha_j)}{(\alpha_i,\alpha_i)}\,,\qquad i,j\in I\, .
 \ee
Note that the Cartan matrix corresponding to $BC_\ell$ is degenerate.
For example the Cartan matrix for $\ell=5$ is given by
\eqref{OSPcar}:
 \be
 A\,= \left(\begin{smallmatrix}
   2&4&-2&0\\
   1&2&-1&0\\
   -1&-2&2&-1\\
   0&0&-1&2
  \end{smallmatrix}\right).
\ee
Using the general formula \eqref{qHam} for $BC_\ell$-root system we
arrive at the following
 \be\label{BCcanG}
  \CH_2^{BC_\ell}\,
  =\,-\sum_{i\in I}\frac{\pr^2}{\pr q_i^2}\,
  +\,\sum_{i=1}^{\ell-1}g_i^2\,e^{q_i-q_{i+1}}\,
  +\,g_\ell^2\,e^{q_{\ell}}\,
 +\,g_{\ell+1}^2\,e^{2q_{\ell}}.
  \ee

It is clear that  specialization of the quadratic Hamiltonian
\eqref{BCcanG} of the  $BC_\ell$-Toda chain to the case of the
coupling constants $g_i^2=1,\,i\in I$ coincides with the quadratic
Hamiltonian $\eqref{BCcan}$ of $\osp(1|2\ell)$-Toda chain. On the
other hand for generic values of $g_i$ in \eqref{BCcanG} it is not
possible by linear changes of variables $q_i$ to transform the
Hamiltonian \eqref{BCcanG}  into the Hamiltonian with $g_i^2=1$.
Thus $\osp(1|2\ell)$-Toda chain realized a special class of
$BC_\ell$-Toda chains.

There is a question on the underlying reason for this
 phenomenon. It is easy to see that the simple
positive roots \eqref{BC2} of $BC_\ell$ and that of $B_{0,\ell}$
\eqref{OPSroot} are closely related. There are however two
differences. First, the short simple root of $\osp(1|2\ell)$ has odd
parity while in $BC_\ell$ root system it is an even root. Second,
while in the case of super Lie algebra $\osp(1|2\ell)$ the
corresponding root system includes the roots $\pm 2\e_\ell$, these
roots are not simple and thus do not enter the expression for the
corresponding Cartan matrix. If however we formally add the root
$2\epsilon_\ell$ to the set of positive simple roots then the
corresponding Cartan matrix constructed according to \eqref{CAM}
precisely coincides with the Cartan matrix of $BC_\ell$ root system.

The fact that in the case of $\osp(1|2\ell)$ the terms of the
Cartan decomposition \eqref{ospNCartan}  corresponding to short
roots are odd actually does not manifest itself in the expressions
 for the Hamiltonians of the corresponding Toda chain.
Indeed, according to Lemma \ref{INDEP} the eigenvalues
$\xi_{\a_i}^{\pm}$ in \eqref{whittvecR}, \eqref{whittvecL} enter the
expressions for quantum  Hamiltonians only via combinations
$g_i^2=\xi_{\a_i}^+\xi_{\a_i}^-$. Therefore  $B_{0,\ell}$-Toda
chains turns out to be a special case of  $BC_\ell$-Toda chain. We
have checked this explicitly for quadratic Hamiltonian in the
previous Section 3.

It is natural to  wonder whether we might  treat the super Lie
algebra $\osp(1|2\ell)$ as a proper  candidate for the Lie algebra
structure associated with $BC_\ell$ root system. Such identification
has at least one  obvious caveat. The root system $BC_\ell$ allows
embedding of roots systems $B_\ell$ and $C_\ell$ having isomorphic
Weyl groups $W_{B_\ell}\simeq W_{C_\ell}$. It would be natural to
expect the same property for the corresponding Lie algebras i.e. a
candidate for the Lie algebra associated with $BC_\ell$ should allow
an embedding of the Lie algebras $\mathfrak{so}_{2\ell+1}$ and
$\mathfrak{sp}_{2\ell}$ associated with the roots systems $B_\ell$
and $C_\ell$ correspondingly. While there indeed exists an embedding
$\mathfrak{sp}_{2\ell}\subset \osp(1|2\ell)$ the super Lie algebra
$\osp(1|2\ell)$ does not allow an embedding of
$\mathfrak{so}_{2\ell+1}$.

\noindent {\small {\bf A.A.G.} {\sl Laboratory for Quantum Field Theory
and Information},\\
\hphantom{xxxx} {\sl Institute for Information
Transmission Problems, RAS, 127994, Moscow, Russia};\\
\hphantom{xxxx} {\it E-mail address}: {\tt anton.a.gerasimov@gmail.com}}\\
\noindent{\small {\bf D.R.L.}
{\sl Laboratory for Quantum Field Theory
and Information},\\
\hphantom{xxxx}  {\sl Institute for Information
Transmission Problems, RAS, 127994, Moscow, Russia};\\
\hphantom{xxxx} {\sl Moscow Center for Continuous Mathematical
Education,\\
\hphantom{xxxx} 119002,  Bol. Vlasyevsky per. 11, Moscow, Russia};\\
\hphantom{xxxx} {\it E-mail address}: {\tt lebedev.dm@gmail.com}}\\
\noindent{\small {\bf S.V.O.} {\sl
 School of Mathematical Sciences, University of Nottingham\,,\\
\hphantom{xxxx} University Park, NG7\, 2RD, Nottingham, United Kingdom};\\
\hphantom{xxxx} {\sl
 Institute for Theoretical and Experimental Physics,
117259, Moscow, Russia};\\
\hphantom{xxxx} {\it E-mail address}: {\tt oblezin@gmail.com}}

\end{document}